\documentclass[a4paper,12pt,intlimits,oneside]{amsart}

\usepackage{enumerate}
\usepackage{latexsym,amssymb}

\newcommand{\comment}[1]{}
\numberwithin{equation}{section}

\newcommand{\bN}{{\mathbb N}}
\newcommand{\bR}{{\mathbb R}}

\newcommand\suth{{  \, |  \, }}
\newcommand\bsuth{{ \, \big |  \, }}

\newcommand{\epf}{ $\Box$\medskip}

\newcounter{rea}
\newcounter{rej}
\newcounter{res}
\newtheorem{thm}{Theorem}[section]
\newtheorem{prop}[thm]{Proposition}

\newtheorem{defn}[thm]{Definition}
\newtheorem{remark}[thm]{Remark}
\begin{document}

\title[]{Pre-dual of Fofana's spaces}


\author[H. G. Feichtinger]{Hans G. Feichtinger}
\address{Faculty of Mathematics, University of Vienna (Autriche),
 Oskar-Morgenstern-Patz 1, 1090 Wien.}
\email{{\tt hans.feichtinger@univie.ac.at}}
\author[J. Feuto]{Justin Feuto}
\address{Laboratoire de Math\'ematiques Fondamentales, UFR Math\'ematiques et Informatique, Universit\'e F\'elix Houphou\"et-Boigny-Cocody, 22 B.P 1194 Abidjan 22. C\^ote d'Ivoire}
\email{{\tt justfeuto@yahoo.fr}}

\thanks{Hans G. Feichtinger is grateful to J. Feuto and
I. Fofana for the opportunity to visit Abidjan and give
a course on Banach Gelfand Triples at the Conference Harmonic Analysis and Applications}
\begin{abstract}

It is the purpose of this paper to give a characterization
of the pre-dual of the spaces introduced by I.~Fofana on the
basis of Wiener amalgam spaces. Those spaces have a specific
dilation behaviour similar to the spaces $L^\alpha(\mathbb R^d)$.
The characterization of the pre-dual will be based on the
idea of minimal invariant spaces (with respect to such a
group of dilation operators).

\end{abstract}
\maketitle

\section{Introduction}
Let $d$ be a fix positive integer and $\mathbb R^d$ the $d$-dimensional Euclidean space, equipped with its Lebesgue measure $dx$. For $1\leq q,p\leq \infty$ the amalgam of $L^q$ and $L^p$ is the space $(L^q,\ell^p)$ of function $f:\mathbb R^d\rightarrow\mathbb C$ which are locally in $L^q$ and such that the sequence $\left\{\left\|f\chi_{I_{k}}\right\|_{q}\right\}_{k\in\mathbb Z^d}$ belongs to $\ell^p(\mathbb Z^d)$, where $I_k=k+\left[0,1\right)^d$.
The map $f\mapsto\left\|f\right\|_p$ denotes the usual norm on the Lebesgue space $L^p(\mathbb R^d)$ on $\mathbb R^d$ while  $\chi_{E}$ stands for the characteristic function of the subset $E$ of $\mathbb R^d$.

Amalgam spaces have being introduced by N. Wiener in 1926 (see \cite{W}), but their systematic study began with the work of  F. Holland \cite{Ho} in 1975. Since then, they have been widely studied (see \cite{BS,FS,Fei1,Fei2,To,FFK} and the references therein).
It is easy to see that the usual Lebesgue space $L^q$  coincides with the  amalgam space $(L^q,\ell^q)$. Also, the Lebesgue space $L^q$ is known to be  invariant under
dilations. In fact, for $\rho>0$, the {\it dilation operator} $St^{(q)}_{\rho}:f\mapsto \rho^{-\frac{d}{q}}f(\rho^{-1}\cdot)$ is isometric. Proper amalgam spaces don't have this property. Even worse, for $q\neq p$ we can't found $\alpha>0$ such that $\sup_{\rho>0}\left\|St^{(\alpha)}_{\rho}f\right\|_{q,p}<\infty$,
although $St^{(\alpha)}_{\rho}f\in (L^q,\ell^p)$ for all $f\in(L^q,\ell^p)$, $\rho>0$ and $\alpha>0$ (see e.g.\ \cite{coni08} or
\cite{ho19}). In order to compensate this shortfall, Fofana introduced (see \cite{Fo1}) in the years 1980's the functions spaces $(L^q,\ell^p)^\alpha$, which consist of $f\in(L^q,\ell^p)$ satisfying $\sup_{\rho>0}\left\|St^{(\alpha)}_{\rho}f\right\|_{q,p}<\infty$ (see Section 2 for more precision).

These spaces can be viewed as some generalized Morrey spaces, and we will always
refer to them as {\it Fofana spaces}.

Many classical results for Lebesgue  and  the classical Morrey spaces have been extended to the setting of the spaces $(L^q,\ell^p)^\alpha(\mathbb R^d)$ (see \cite{Fo1,Fo2,Fo3,Fe1,Fe2,Fe3,KF,DFS,SF,CF}).
Although the dual space of Lebesgue spaces $L^\alpha(\mathbb R^d)$ ($1\leq \alpha<\infty)$ and amalgam spaces $(L^q,\ell^p)(\mathbb R^d)$ ($1\leq q,p<\infty$) are well known ($L^{\alpha'}$ and $(L^{q'},\ell^{p'})$ respectively with $\frac{1}{p'}+\frac{1}{p}=1$), the one of $(L^q,\ell^p)^\alpha(\mathbb R^d)$ is still unknown. But in the case $p=\infty$ and $q<\alpha$, there are already four characterizations of the pre-dual of $(L^q,\ell^\infty)^\alpha(\mathbb R^d)$ (see \cite{AX,GM,K,Z}) which are equal.

The purpose of this paper is to give a beginning of answer by the determination of a pre-dual of $(L^q,\ell^p)^\alpha(\mathbb R^d)$ for $1<q<\alpha<p<\infty$.

In doing so we will make use of the idea of minimal invariant Banach spaces of functions which has already quite some tradition and has shown to be useful in a variety of situations, see \cite{FEI77-3,FEI81,FEI81-2,FEI87-1,FEI88} or \cite{JO79}.

The paper is organized as follows :
In Section 2, we give some basic facts about {\it amalgams} and minimal Banach spaces.
The third section is devoted to  pre-dual of {\it Fofana spaces}
as well as certain properties of these spaces.

\section{Some basic facts about Amalgam and Minimal Banach spaces}

For any normed  space $E$, we denote by $E^\ast$ its topological dual space. 
Given $1\leq p,q\leq\infty$, the amalgam space $(L^q,\ell^p)$ is equipped with the norm $\left\|f\right\|_{q,p}=\left\|\left\{\left\|f\chi_{I_k}\right\|_{q}\right\}\right\|_{\ell^p}$.
For any $\rho>0$, we put
\begin{equation}
_{\rho}\left\|f\right\|_{q,p}=\left\{\begin{array}{lll}\left[\sum_{k\in\mathbb Z^d}\left\|f\chi_{I^{\rho}_{k}}\right\|^{p}_{q}\right]^{\frac{1}{p}}&\text{ if }&p<\infty\\
\sup_{k\in\mathbb Z^d}\left\|f\chi_{I^{\rho}_{k}}\right\|_q&\text{ if }&p=\infty\end{array}\right.\label{1.1}
\end{equation}
with
$I^{\rho}_{k}=\Pi^{d}_{j=1}\left[k_j\rho,(k_j+1)\rho\right)\text{ if }  k=(k_j)_{1\leq j\leq d}\in\mathbb Z^{d}.$  It is clear that $\left\|f\right\|_{q,p}=\ _{1}\left\|f\right\|_{q,p}$.

We have the following well known properties (see for example \cite{FS}).
\begin{enumerate}
\item For $0<\rho<\infty$, $f\mapsto\ _{\rho}\left\|f\right\|_{q,p}$ is a norm on $(L^q,\ell^p)(\mathbb R^d)$ equivalent to $f\mapsto\left\|f\right\|_{q,p}$.
     With respect to these norms the
    amalgam spaces  $(L^q,\ell^p)(\mathbb R^d)$ are Banach spaces.
\item The spaces are (strictly) increasing with the global exponent $q$ and (strictly) decreasing with a growing local exponent $q$; more precisely
\begin{eqnarray}
\left\|f\right\|_{q,p}\leq\left\|f\right\|_{q_1,p}&\text{ if }&q<q_1\leq\infty,\label{1.4}\\
\left\|f\right\|_{q,p}\leq\left\|f\right\|_{q,p_1}&\text{ if }&1\leq p_1<p.\label{1.5}
\end{eqnarray}
\item For $0<\rho<\infty$, H\"older's inequality is fulfilled :
\begin{equation}
\left\|fg\right\|_1\leq\ _{\rho}\left\|f\right\|_{q,p}\ _{\rho}\left\|g\right\|_{q',p'},\ \ f,g\in L^{1}_{loc}(\mathbb R^d),\label{1.6}
\end{equation}
where $q'$ and $p'$ are conjugate exponent of $q$ and $p$ respectively: $\frac{1}{q}+\frac{1}{q'}=1=\frac{1}{p}+\frac{1}{p'}$. When $1\leq q,p<\infty$, $(L^{q'},\ell^{p'})(\mathbb R^d)$ is isometrically isomorphic to the dual $(L^{q},\ell^p)(\mathbb R^d)^{\ast}$ of $(L^{q},\ell^p)(\mathbb R^d)$ in the sense that for any element $T$ of $(L^{q},\ell^p)(\mathbb R^d)^{\ast}$, there is an unique element $\phi(T)$ of $(L^{q'},\ell^{p'})(\mathbb R^d)$ such that
$$\left\langle T,f\right\rangle=\int_{\mathbb R^d}f(x)\phi(T)(x)dx,\ \ f\in(L^{q},\ell^p)(\mathbb R^d)$$
and furthermore
\begin{equation}
\left\|\phi(T)\right\|_{q',p'}=\left\|T\right\|.\label{1.7}
\end{equation}
We recall that $\left\|T\right\|:=\sup\left\{\left|\left\langle T,f\right\rangle\right|
 \bsuth f\in L^{q}_{loc}(\mathbb R^d)\text{ and } \left\|f\right\|_{q,p}\leq 1\right\}.$
\end{enumerate}

Next we summarize a couple of properties of dilation operators.
We assume that $1\leq\alpha\leq\infty$.
\begin{enumerate}
\item For any real number $\rho>0$, $St^{(\alpha)}_{\rho}$ maps $L^{1}_{loc}(\mathbb R^d)$ into itself.
\item $St^{(\alpha)}_{1}f=f,\ f\in L^{1}_{loc}(\mathbb R^d)$.
\item $St^{(\alpha)}_{\rho_1}\circ St^{(\alpha)}_{\rho_2}=St^{(\alpha)}_{\rho_1\rho_2}=St^{(\alpha)}_{\rho_2}\circ St^{(\alpha)}_{\rho_1}$.
\end{enumerate}
In other words, $(St^{(\alpha)}_{\rho})_{\rho>0}$ is a commutative group of operators on $L^{1}_{loc}(\mathbb R^d)$, isomorphic to the multiplicative group $(0,\infty)$. As mentioned in the introduction, we have for $1\leq \alpha\leq \infty$
\begin{equation}\left\|St^{(\alpha)}_{\rho}f\right\|_\alpha=\left\|f\right\|_\alpha, \ 0<\rho<\infty\text{ and }f\in L^{1}_{loc}(\mathbb R^d). \,\, \label{2.4}
\end{equation}
In other words, each of those normalizations is isometric on exactly one of the family of $L^r$-spaces.

For amalgam spaces, direct computations (see for example (2.1) and Proposition 2.2 in \cite{Fe3})
give the following results:
\begin{equation*}
\left\|St^{(\alpha)}_{\rho}f\right\|_{q,p}=\rho^{-d(\frac{1}{\alpha}-\frac{1}{q})}\ _{\rho^{-1}}
\negthinspace \left\|f\right\|_{q,p},\ 0<\rho<\infty\text{ and }f\in L^{1}_{loc}(\mathbb R^d)
\end{equation*}
and therefore
\begin{equation}
\left\|f\right\|_{q,p,\alpha}:=\sup_{\rho>0}\left\|St^{(\alpha)}_{\rho}f\right\|_{q,p}=\sup_{\rho>0}\rho^{d(\frac{1}{\alpha}-\frac{1}{q})}\ _{\rho}\left\|f\right\|_{q,p}.\label{2.6}
\end{equation}
It follows that the space $(L^q,\ell^p)^\alpha(\mathbb R^d)$ can be defined by
\begin{equation}
(L^q,\ell^p)^\alpha(\mathbb R^d)=\left\{f\in L^{1}_{loc}(\mathbb R^d) \bsuth \left\|f\right\|_{q,p,\alpha}<\infty\right\}.
\end{equation}
We recall that for $q<\alpha$, the space $(L^q,\ell^\infty)^\alpha(\mathbb R^d)$ is exactly the classical Morrey space $\mathcal M_{q,\frac{dq}{\alpha}}(\mathbb R^d)$ introduced by Morrey in 1938, see \cite{M}, and defined for $\lambda=\frac{dq}{\alpha}$ by
\begin{equation}
\mathcal M_{q,\lambda}(\mathbb R^d)=\left\{f\in L^{q}_{loc}(\mathbb R^d)  \bsuth
\sup_{x\in\mathbb R^d,r>0}r^{\frac{\lambda-d}{q}}\left\|f\chi_{B(x,r)}\right\|_{q}<\infty\right\}.
\end{equation}
Fofana spaces have the following properties (cf.  \cite{Fo1},  \cite{Fo2}):
\begin{enumerate}
\item $\left((L^q,\ell^p)^\alpha(\mathbb R^d),\left\|\cdot\right\|_{q,p,\alpha}\right)$ is a Banach space which is non trivial if and only if $q\leq\alpha\leq p$,
\item if $\alpha\in \left\{p,q\right\}$ then $(L^q,\ell^p)^\alpha(\mathbb R^d)=L^\alpha(\mathbb R^d)$  with equivalent norms,
\item if $q<\alpha<p$ then $L^\alpha(\mathbb R^d)\subsetneq L^{\alpha,\infty}(\mathbb R^d)\subsetneq(L^q,\ell^p)^\alpha(\mathbb R^d)\subsetneq (L^q,\ell^p)(\mathbb R^d)$, where $L^{\alpha,\infty}(\mathbb R^d)$ is the weak Lebesgue space on $\mathbb R^d$ defined by
$$L^{\alpha,\infty}(\mathbb R^d)=\left\{f\in L^{1}_{loc}(\mathbb R^d) \bsuth \left\|f\right\|^{\ast}_{\alpha,\infty}<\infty\right\},$$
with $\left\|f\right\|^{\ast}_{\alpha,\infty}:=\sup_{\lambda>0}\left|\left\{x\in\mathbb R^d \suth \left|f(x)\right|>\lambda\right\}\right|^{\frac{1}{\alpha}}$. We denote by  $\left|E\right|$, the Lebesgue measure of a measurable subset $E$ of $\mathbb R^d$.
\end{enumerate}
For any $\rho>0$, the dilation $St^{(\alpha)}_{\rho}$ map isometrically the space $(L^q,\ell^p)^\alpha(\mathbb R^d)$ to itself. More precisely,
\begin{equation}
\left\|St^{(\alpha)}_{\rho}f\right\|_{q,p,\alpha}=\left\|f\right\|_{q,p,\alpha},\ 0<\rho<\infty\text{ and }f\in L^{1}_{loc}(\mathbb R^d).
\end{equation}
\begin{remark}
 It is easy to see that, if $0<\rho<\infty$ and $(f,g)$ is an element of $L^{1}_{loc}(\mathbb R^d)\times L^{1}_{loc}(\mathbb R^d)$ such that $(St^{(\alpha)}_{\rho}f)g$ belongs to $L^1(\mathbb R^d)$, then $f(St^{(\alpha)}_{\rho^{-1}}g)$ belongs to $L^1(\mathbb R^d)$ and
\begin{equation}
\int_{\mathbb R^d}\left(St^{(\alpha)}_{\rho}f\right)(x)g(x)dx=\int_{\mathbb R^d}f(x)\left(St^{(\alpha')}_{\rho^{-1}}g\right)(x)dx.\label{2.8}
\end{equation}
\end{remark}
Searching for a dilation invariant version of the Segal algebra $S_0(\bR^d)$
 (introduced in \cite{FEI81-2}) Feichtinger and Zimmermann introduced a certain
 {\it exotic minimal space} in \cite{FZ}.
We will use later on 
the following result of that paper:
\begin{thm}[\cite{FZ}, Theorem 2.1] \label{minbanach}Let $\bf B$ be a Banach space, and $\Phi=(\varphi_j)_{j\in J}$ a (not necessarily countable) bounded family in $\bf B$. Define
$$\mathcal B=\mathcal B_\Phi:=\left\{f=\sum_{j\in J}a_j\varphi_j:\sum_{j\in J}\left|a_j\right|<\infty\right\},$$
and let
$$\left\|f\right\|_{\mathcal B}=\inf\left\{\sum_{j\in J}\left|a_j\right| \bsuth
f=\sum_{j\in J}a_j\varphi_j\right\}.$$
Then $(\mathcal B,\left\|\cdot\right\|_{\mathcal B})$ is a Banach space continuously embedded into $\bf B$.
\end{thm}

\vspace{2mm}

\section{A pre-dual of Fofana spaces}

\begin{defn}
Let $1\leq q\leq \alpha\leq p\leq \infty$. The space $\mathcal H(q,p,\alpha)$ is defined as the set of all elements $f$ of $L^{1}_{loc}(\mathbb R^d)$ for which there exist a sequence $\left\{(c_n,\rho_n,f_n)\right\}_{n\geq 1}$ of elements of $\mathbb C\times (0,\infty)\times(L^{q'},\ell^{p'})(\mathbb R^d)$ such that
\begin{equation}\left\{\begin{array}{l}\sum_{n\geq 1}\left|c_n\right|<\infty\\
\left\|f_n\right\|_{q',p'}\leq 1,\  n\geq 1\\
f:=\sum_{n\geq 1}c_nSt^{(\alpha')}_{\rho_n}f_n\text{ in the sense of }L^1_{loc}(\mathbb R^d)\end{array}\right..\label{decomposition}
\end{equation}
\end{defn}
We will always refer to any sequence $\left\{(c_n,\rho_n,f_n)\right\}_{n\geq 1}$ of elements of $\mathbb C\times (0,\infty)\times(L^{q'},\ell^{p'})(\mathbb R^d)$ satisfying (\ref{decomposition}) as $\mathfrak h$-decomposition of $f$.

For any element $f$ of $\mathcal H(q,p,\alpha)$, we set
\begin{equation}
\left\|f\right\|_{\mathcal H(q,p,\alpha)}:=\inf\left\{\sum_{n\geq 1}\left|c_n\right|\right\},\label{2.9}
\end{equation}
where the infimum is taken over all $\mathfrak h$-decomposition of $f$.
\begin{prop}
Let $1\leq q\leq \alpha\leq p\leq \infty$.   $\mathcal H(q,p,\alpha)$ endowed with $\left\|\cdot\right\|_{\mathcal H(q,p,\alpha)}$ is a Banach space continuously embedded into $L^{\alpha'}(\mathbb R^d)$.
\end{prop}
\proof Since $1\leq p'\leq \alpha'\leq q'\leq\infty$, it comes from (\ref{1.4}) and (\ref{1.5}) that for any element $f\in (L^{q'},\ell^{p'})(\mathbb R^d)$,
$$f\in L^{\alpha'}\text{ and }\left\|f\right\|_{\alpha'}\leq\left\|f\right\|_{q',p'}.$$
Therefore, by (\ref{2.4})
$$St^{(\alpha')}_{\rho}f\in L^{\alpha'} \text{ and }\left\|St^{(\alpha')}_{\rho}f\right\|_{\alpha'}\leq \left\|f\right\|_{q',p'}.$$
The result follows from Theorem \ref{minbanach}, using the Banach space
${\bf B}=L^{\alpha'}(\mathbb R^d)$ and its bounded subset
$$\Phi=\left\{St^{(\alpha')}_{\rho}f \bsuth \rho>0 \text{ and }\left\|f\right\|_{q',p'}\leq 1\right\}.$$
\epf
\medskip

For $1<p<\infty$ and $0<\lambda<d$. Zorko proved in \cite{Z} that the Morrey space $\mathcal M_{p,\lambda}$ is the dual space of the space $\mathcal Z_{p',\lambda}$, which consists of all the functions $f$ on $\mathbb R^d$ which can be written in the form $f:=\sum_k c_k \mathfrak{a}_k$, where $\left\{c_k\right\}$ is a sequence in $\ell^1$, and $\left\{\mathfrak{a}_k\right\}$ is a sequence of functions on $\mathbb R^d$ satisfying for each $k$,
\begin{itemize}
    \item $\mathrm{supp } \ \mathfrak{a}_k\subset \text{ a ball } B_k$,
    \item $ \left\|\mathfrak{a}_k\right\|_{p'}\leq
        1/{\left|B_k\right|^{\frac{d-\lambda}{dp}}}.$
    \end{itemize}
		Notice that for $p=\infty$ and $q<\alpha$, the space $\mathcal Z_{q',\frac{dq}{\alpha}}$ is continuously embedded into $\mathcal H(q,\infty,\alpha)$. In fact, let $f=\sum_k c_k \mathfrak{a}_k\in\mathcal Z_{q',\frac{dq}{\alpha}}$ with the sequences $\left\{c_k\right\}$  and $\left\{\mathfrak{a}_k\right\}$ as in the Introduction. For  any   $k \in \bN$, we put
$$u_k=St^{(\alpha')}_{r^{-1}_k}(\mathfrak{a}_k),$$
where $r_k$ is the radius of the ball $B_k$ associated to $\mathfrak{a}_k$. There exists a constant $C>0$ (one can take $C=\left|B(0,1)\right|^{\frac{1}{\alpha}-\frac{1}{q}}$), such that
$$\left\|2^{-d}C^{-1}u_k\right\|_{q',1}\leq 1\text{ and }f=\sum_{k\geq 1}(2^dCc_k)St^{(\alpha')}_{r_k}(2^{-d}C^{-1}u_k).$$

We prove   next  that for $1\leq q\leq \alpha\leq p\leq\infty$, the space $\mathcal H(q,p,\alpha)$ is a certain minimal Banach space.
\begin{prop}
For $1\leq q\leq\alpha\leq p\leq\infty$, the space $\mathcal H(q,p,\alpha)$ is a minimal Banach space, isometrically invariant under the family $(St^{(\alpha')}_{\rho})_{\rho>0}$ and such that
$$(L^{q'},\ell^{p'})(\mathbb R^d)\subset \mathcal H(q,p,\alpha)\subset L^{1}_{loc}(\mathbb R^d),$$  with continuous inclusions.
 \end{prop}
\proof\begin{enumerate}\item Let us first prove that $\mathcal H(q,p,\alpha)$ is isometrically invariant by $St^{(\alpha')}_{\rho}$ for all $\rho>0$.

Let $f\in \mathcal H(q,p,\alpha)$. For $0<\rho<\infty$ and any $\mathfrak h$-decomposition $\left\{(c_n,\rho_n,f_n)\right\}_{n\geq 1}$ of $f$, we have
$$St^{(\alpha')}_{\rho}f=\sum_{n\geq 1}c_nSt^{(\alpha')}_{\rho\rho_n}f_n$$
so that
\begin{equation}
St^{(\alpha')}_{\rho}f\in \mathcal H(q,p,\alpha)\text{ and }\left\|St^{(\alpha')}_{\rho}f\right\|_{\mathcal H(q,p,\alpha)}=\left\|f\right\|_{\mathcal H(q,p,\alpha)}.\label{2.10}
\end{equation}
Thus, $St^{(\alpha')}_{\rho}$ is an isometric automorphism of $\mathcal H(q,p,\alpha)$.
\item First we verify that $(L^{q'},\ell^{p'})$ is continuously embedded into $\mathcal H(q,p,\alpha)$.
For any   $0 \neq f \in(L^{q'},\ell^{p'})$ we have
$$f=\left\|f\right\|_{q',p'}St^{(\alpha_1)}_{1}(\left\|f\right\|^{-1}_{q',p'}f)\text{  and }\left\|\left\|f\right\|^{-1}_{q',p'}f\right\|_{q',p'}=1$$
and therefore $f$ belongs to $\mathcal H(q,p,\alpha)$ and satisfies
\begin{equation}
\left\|f\right\|_{\mathcal H(q,p,\alpha)}\leq\left\|f\right\|_{q',p'}.\label{2.11}
\end{equation}
Thus our claim is verified.
\item It remains to prove that this space is minimal.
 Let $X$ be another Banach space continuously embedded into $L^{1}_{loc}(\mathbb R^d)$ such that
\begin{enumerate}
\item $X$ is isometrically invariant by $(St^{(\alpha')}_{\rho})_{\rho>0}$, i.e., if $f\in X$
then $St^{(\alpha')}_{\rho}f\in X\text{ and }\left\|St^{(\alpha')}_{\rho}f\right\|_X=\left\|f\right\|_X,\ \rho>0$,
\item $(L^{q'},\ell^{p'})(\mathbb R^d)$ is continuously embedded into $X$, i.e., there exists $K>0$ such that for $f\in (L^{q'},\ell^{p'})(\mathbb R^d)$, $f\in X$ and $\left\|f\right\|_{X}\leq K\ _{1}\left\|f\right\|_{q',p'}$.
\end{enumerate}
 For  $f \in\mathcal H(q,p,\alpha)$ and any $\mathfrak h$-decomposition $\left\{(c_n,\rho_n,f_n)\right\}_{n\geq 1}$  
$$c_nSt^{(\alpha')}_{\rho_n}f_n\in X\text{ and }\left\|c_nSt^{(\alpha')}_{\rho_n}f_n\right\|_{X}=\left|c_n\right|\left\|f_n\right\|_X\leq K\left|c_n\right|,\ n\geq1.$$
Thus
$$\sum_{n\geq 1}\left\|c_nSt^{(\alpha')}_{\rho_n}f_n\right\|_{X}\leq K\sum_{n\geq 1}\left|c_n\right|$$
and therefore, as $X$ is a Banach space,
$$f=\sum_{n\geq 1}c_n St^{(\alpha')}_{\rho_n}f_n\in X\text{ and }\left\|f\right\|_{X}\leq K\sum_{n\geq 1}\left|c_n\right|.$$
From this and (\ref{2.9}) it follows that
$$f\in X\text{ and }\left\|f\right\|_{X}\leq K\left\|f\right\|_{\mathcal H(q,p,\alpha)}.$$
Thus $\mathcal H(q,p,\alpha)$ is continuously included in $X$.
\end{enumerate}
\epf

Our main result can be stated as follows.
\begin{thm}\label{main}
Let $1<q\leq\alpha\leq p\leq\infty$. The operator $g\mapsto T_g$ defined by (\ref{2.13}) is an isometric isomorphism of $(L^q,\ell^p)^\alpha(\mathbb R^d)$ into $\mathcal H(q,p,\alpha)^{\ast}$.
\end{thm}
For the proof, we need some intermediate results.
\begin{remark}\label{rem2.3}
It is clear that if $\left\{(c_n,\rho_n,f_n)\right\}_{n\geq 1}$ is 
a $\mathfrak h$-decomposition of $f \in \mathcal H(q,p,\alpha)$,
then $(\sum^{m}_{n= 1}c_n St^{(\alpha')}_{\rho_n}f_n)_{m\geq 1}$ is a sequence of elements of $(L^{q'},\ell^{p'})(\mathbb R^d)$ which converges to $f$ in $\mathcal H(q,p,\alpha)$. Hence $(L^{q'},\ell^{p'})(\mathbb R^d)$ is a dense subspace of $\mathcal H(q,p,\alpha)$.\label{b}
\end{remark}

\begin{prop}\label{prop2.4}
Let $1\leq q\leq\alpha\leq p\leq \infty$, $f\in \mathcal H(q,p,\alpha)$ and $g\in (L^q,\ell^p)^\alpha$. Then $fg$ belongs to $L^1(\mathbb R^d)$ and
\begin{equation}
\left|\int_{\mathbb R^d}f(x)g(x)dx\right|\leq \left\|f\right\|_{\mathcal H(q,p,\alpha)}\left\|g\right\|_{q,p,\alpha}.\label{2.12}
\end{equation}
\end{prop}
\proof
Let $\left\{(c_n,\rho_n,f_n)\right\}_{n\geq 1}$ be a $\mathfrak h$-decomposition of $f$.

By using successively (\ref{2.8}), (\ref{1.6}) and (\ref{2.6}), we obtain for any  $n\geq1$,
\begin{eqnarray*}
\left|\int_{\mathbb R^d}St^{(\alpha')}_{\rho_n}f_n(x)g(x)dx\right|&=&\left|\int_{\mathbb R^d}f_n(x)St^{(\alpha)}_{\rho^{-1}_n}g(x)dx\right|\leq\int_{\mathbb R^d}\left|f_n(x)St^{(\alpha)}_{\rho^{-1}_n}g(x)\right|dx\\
&\leq&\left\|f_n\right\|_{q',p'}\left\|St^{(\alpha)}_{\rho^{-1}_{n}}g\right\|_{q,p}\leq \left\|St^{(\alpha)}_{\rho^{-1}_{n}}g\right\|_{q,p}\leq \left\|g\right\|_{q,p,\alpha}.
\end{eqnarray*}
Therefore we have
$$\sum_{n\geq 1}\int_{\mathbb R^d}\left|c_nSt^{(\alpha')}_{\rho_n}f_n(x)g(x)\right|dx\leq \left\|g\right\|_{q,p,\alpha}\sum_{n\geq 1}\left|c_n\right|.$$
This implies that $fg=\sum_{n\geq 1}c_nSt^{(\alpha')}_{\rho_n}f_ng$ belongs to $L^1(\mathbb R^d)$ and
$$\left|\int_{\mathbb R^d}f(x)g(x)dx\right|\leq \int_{\mathbb R^d}\left|f(x)g(x)\right|dx\leq \left\|g\right\|_{q,p,\alpha}\sum_{n\geq 1}\left|c_n\right|.$$
Taking the infimum with respect to all $\mathfrak h$-decompositions of $f$,
we get
$$\left|\int_{\mathbb R^d}f(x)g(x)dx\right|\leq \int_{\mathbb R^d}\left|f(x)g(x)\right|dx\leq \left\|g\right\|_{q,p,\alpha}\left\|f\right\|_{H(q,p,\alpha)}$$
\epf

\begin{remark}\label{rem2.5}
Let $1\leq q\leq \alpha\leq p\leq \infty$ and set :
\begin{equation}
\left\langle T_g,f\right\rangle=\int_{\mathbb R^d}f(x)g(x)dx,\ g\in(L^q,\ell^p)^\alpha(\mathbb R^d)\text{ and }f\in \mathcal H(q,p,\alpha).\label{2.13}
\end{equation}
By Prop. \ref{prop2.4} and the fact that $\varphi\mapsto\int_{\mathbb R^d}\varphi(x)dx$ belongs to $L^{1}(\mathbb R^d)^{\ast}$, that :
\begin{enumerate}
\item for any element $g$ of $(L^q,\ell^p)^\alpha(\mathbb R^d)$, $f\mapsto\left\langle T_g,f\right\rangle$ belongs to $\mathcal H(q,p,\alpha)^{\ast}$,
\item $T:g\mapsto T_g$ is linear and bounded mapping from  $(L^q,\ell^p)^\alpha(\mathbb R^d)$ into $\mathcal H(q,p,\alpha)^{\ast}$ satisfying $\left\|T\right\|\leq 1$, that is :
\begin{equation}
\left\|T_g\right\|\leq\left\|g\right\|_{q,p,\alpha},\ g\in(L^q,\ell^p)^\alpha(\mathbb R^d).\label{2.14}
\end{equation}
\end{enumerate}
\end{remark}
Now we can prove  our main result.
\proof[Proof of Theorem \ref{main}]
We know (see Remark \ref{rem2.5}) that $g\mapsto T_g$ is a bounded linear application of $(L^q,\ell^p)^\alpha(\mathbb R^d)$ into $\mathcal H(q,p,\alpha)^{\ast}$ such that
$$\left\|T_g\right\|\leq \left\|g\right\|_{q,p,\alpha},\ g\in(L^q,\ell^p)^\alpha(\mathbb R^d).$$
\label{1}
 Let $T$ be an element of $\mathcal H(q,p,\alpha)^{\ast}$. From (\ref{2.11}) it follows that the restriction $T_0$ of $T$ to $(L^{q'},\ell^{p'})(\mathbb R^d)$ belongs to $(L^{q'},\ell^{p'})(\mathbb R^d)^{\ast}$. Furthermore, we have $1\leq p'\leq \alpha'\leq q'<\infty$. So, by (\ref{1.7}), there is an element $g$ of $(L^q,\ell^p)(\mathbb R^d)$ such that
\begin{equation}
\left\langle T,f\right\rangle =\int_{\mathbb R^d}f(x)g(x)dx,\ f\in(L^{q'},\ell^{p'})(\mathbb R^d).\label{**}
\end{equation}
Hence , for   $f \in (L^{q'},\ell^{p'})(\mathbb R^d)$ and  $\rho>0$ we have
$$\int_{\mathbb R^d}St^{(\alpha)}_{\rho}g(x)f(x)dx=\int_{\mathbb R^d}g(x)St^{(\alpha')}_{\rho^{-1}}f(x)dx=\left\langle T,St^{(\alpha')}_{\rho^{-1}}f\right\rangle$$
by (\ref{1.6}) and (\ref{2.8}), and
$$\left|\int_{\mathbb R^d}St^{(\alpha)}_{\rho}g(x)f(x)dx\right|\leq \left\|T\right\|\left\|St^{(\alpha')}_{\rho^{-1}}f\right\|_{\mathcal H(q,p,\alpha)}=\left\|T\right\|\left\|f\right\|_{\mathcal H(q,p,\alpha)}\leq \left\|T\right\|\left\|f\right\|_{q',p'}$$
by (\ref{2.10}) and (\ref{2.11}). From this and (\ref{1.7}) it follows that
$$St^{(\alpha)}_{\rho}g\in (L^q,\ell^p)(\mathbb R^d)\text{ and }\left\|St^{(\alpha)}_{\rho}g\right\|_{q,p}\leq\left\|T\right\|,\ \rho\in(0,\infty)$$
and therefore, by (\ref{2.6}),
$$\left\|g\right\|_{q,p,\alpha}\leq \left\|T\right\|\text{ and }g\in(L^q,\ell^p)^\alpha(\mathbb R^d).$$
From (\ref{**}), the density of $(L^{q'},\ell^{p'})(\mathbb R^d)$ in $\mathcal H(q,p,\alpha)$ (see Remark \ref{rem2.3}) and Prop. \ref{prop2.4}, we get
$$\left\langle T,f\right\rangle=\int_{\mathbb R^d}f(x)g(x)dx,\ f\in \mathcal H(q,p,\alpha).$$
This completes the proof.
\epf

We end this paper by stating  some interesting properties of the spaces $\mathcal H(q,p,\alpha)$.
\begin{prop}Let $1\leq q\leq\alpha\leq p\leq \infty$.
\begin{enumerate}
\item The space $\mathcal H(q,p,\alpha)$ is a Banach $L^1(\mathbb R^d)$-module: \\ for $(f,\varphi)\in \mathcal H(q,p,\alpha)\times L^1(\mathbb R^d)$
\begin{equation}
f\ast\varphi\in \mathcal H(q,p,\alpha)\text{ and }\left\|f\ast\varphi\right\|_{\mathcal H(q,p,\alpha)}\leq \left\|f\right\|_{\mathcal H(q,p,\alpha)}\left\|\varphi\right\|_1.\label{2.15}
\end{equation}
\item If $1<q$ and $\varphi\in L^1(\mathbb R^d)$ with $\varphi\geq 0$ and $\left\|\varphi\right\|_1=1$ then
\begin{equation}
\lim_{\epsilon\rightarrow 0}\left\|f\ast St^{(1)}_{\epsilon}\varphi-f\right\|_{\mathcal H(q,p,\alpha)}=0,\ \ f\in \mathcal H(q,p,\alpha). \end{equation}
Consequently $\mathcal H(q,p,\alpha)$ is an essential  Banach $L^1(\mathbb R^d)$-module.
\item  For $1<p,q \leq \infty$ the Schwartz space
of rapidly decreasing functions or the space of compactly supported $C^\infty$-functions are
dense subspaces of $ \mathcal H(q,p,\alpha)$.
\item  For $1<p,q \leq \infty$, the Banach space $ \mathcal H(q,p,\alpha)$ is separable.
\end{enumerate}
\end{prop}
\proof
Let $f\in\mathcal H(q,p,\alpha)$ and $\varphi\in L^1(\mathbb R^d)\setminus\left\{0\right\}$ .
\begin{enumerate}
\item 
Let $\left\{(c_n,\rho_n,f_n)\right\}_{n\geq 1}$ be a $\mathfrak h$-decomposition of $f$.
A direct computation shows that
$$(St^{(\alpha')}_{\rho_n}f_n)\ast\varphi=St^{(\alpha')}_{\rho_n}\left[f_n\ast St^{(1)}_{\rho^{-1}_{n}}\varphi\right],\ n\geq 1.$$
Hence, combining the fact that $(L^{q'},\ell^{p'})(\mathbb R^d)$ is a $L^1(\mathbb R^d)$ Banach module and Relation (\ref{2.4}), we obtain
$$\left\|f_n\ast St^{(1)}_{\rho^{-1}_{n}}\varphi\right\|_{q',p'}\leq \left\|f_n\right\|_{q',p'}\left\|\varphi\right\|_1,\ n\geq 1.$$
This implies that
$$\left\|\left\|\varphi\right\|^{-1}_1f_n\ast St^{(1)}_{\rho^{-1}_{n}}\varphi\right\|_{q',p'}\leq\left\|f_n\right\|_{q',p'}\leq 1,\  n\geq 1.$$
Moreover
\begin{eqnarray*}\sum_{n\geq1}\left\|c_n(St^{(\alpha')}_{\rho_n}f_n)\ast\varphi\right\|_{\alpha'}&\leq&\sum_{n\geq 1}\left|c_n\right|\left\|St^{(\alpha')}_{\rho_n}f_n\right\|_{\alpha'}\left\|\varphi\right\|_1=\sum_{n\geq 1}\left|c_n\right|\left\|f_n\right\|_{\alpha'}\left\|\varphi\right\|_1\\
&\leq&\sum_{n\geq 1}\left|c_n\right|\left\|f_n\right\|_{q',p'}\left\|\varphi\right\|_1\leq(\sum_{n\geq 1}\left|c_n\right|)\left\|\varphi\right\|_1<\infty.
\end{eqnarray*}
Therefore
$$ f\ast\varphi=\sum_{n\geq 1}c_n(St^{(\alpha')}_{\rho_n}f_n)\ast\varphi=\sum_{n\geq 1}c_n\left\|\varphi\right\|_{1}St^{(\alpha')}_{\rho_n}(\left\|\varphi\right\|^{-1}_{1}f_n\ast St^{(1)}_{\rho^{-1}_{n}}\varphi)$$
in the sense of $L^{1}_{loc}(\mathbb R^d)$, with
$$\sum_{n\geq 1}\left|c_n\left\|\varphi\right\|_1\right|=(\sum_{n\geq 1}\left|c_n\right|)\left\|\varphi\right\|_1<\infty\text{ and } \left\|\left\|\varphi\right\|^{-1}_{1}f_n\ast St^{(1)}_{\rho^{-1}_{n}}\varphi\right\|_{q',p'}\leq 1.$$
That is $f\ast\varphi$ belongs to $\mathcal H(q,p,\alpha)$ and satisfies
$$\left\|f\ast\varphi\right\|_{\mathcal H(q,p,\alpha}\leq\left\|f\right\|_{\mathcal H(q,p,\alpha)}\left\|\varphi\right\|_1.$$
\item Let us assume that $1<q$ 
and $\left\|\varphi\right\|_1
=1$. For any real number $\epsilon>0$, we set $\varphi_\epsilon=St^{(1)}_{\epsilon}\varphi$. Let us consider $\left\{(c_n,\rho_n,f_n)\right\}_{n\geq 1} $, a $\mathfrak h$-decomposition of $f$. b
We know that the sequence $(f^{m})_{m\geq 1}$ defined by
$$f^{m}=\sum^{m}_{n=1}c_nSt^{(\alpha')}_{\rho_n}f_n,\ \ m\geq1$$
converges to $f$ in $\mathcal H(q,p,\alpha)$. Let us fix any real number $\delta>0$. There is a positive integer $m_\delta$ satisfying $$\left\|f-f^{m_\delta}\right\|<\frac{\delta}{3}.$$
Moreover, for any real number $\epsilon>0$, we have
\begin{eqnarray*}
\left\|f\ast\varphi_\epsilon-f\right\|_{\mathcal H(q,p,\alpha)}&\leq&\left\|f\ast\varphi_\epsilon-f^{m_\delta}\ast\varphi_\epsilon\right\|_{\mathcal H(q,p,\alpha)}+\left\|f^{m_\delta}\ast\varphi_\epsilon-f^{m_\delta}\right\|_{\mathcal H(q,p,\alpha)}\\
&+&\left\|f^{m_\delta}-f\right\|_{\mathcal H(q,p,\alpha)}\\
&\leq&2\left\|f^{m_\delta}-f\right\|_{\mathcal H(q,p,\alpha)}+\left\|\sum^{m_\delta}_{n=1}c_n(St^{(\alpha')}_{\rho_n}f_n)\ast\varphi_\epsilon-\sum^{m_\delta}_{n=1}c_nSt^{(\alpha')}_{\rho_n}f_n\right\|_{\mathcal H(q,p,\alpha)}.
\end{eqnarray*}
The last term is not greater than
$$\frac{2\delta}{3}+\left\|\sum_{n\in N(m_\delta)}c_n\left\|f_n\ast\varphi_{\epsilon\rho^{-1}_{n}}-f_n\right\|_{q',p'}St^{(\alpha')}_{\rho_n}\left(\frac{f_n\ast\varphi_{\epsilon\rho^{-1}_{n}}-f_n}{\left\|f_n\ast\varphi_{\epsilon\rho^{-1}_{n}}-f_n\right\|_{q',p'}}\right)\right\|_{\mathcal H(q,p,\alpha)},$$
where
$$N(m_\delta)=\left\{n/1\leq n\leq m_\delta\text{ and }\left\|f_n\ast\varphi_{\epsilon\rho^{-1}_{n}}-f_n\right\|_{q',p'}\neq 0\right\}.$$
It comes that
$$\left\|f\ast\varphi_\epsilon-f\right\|_{\mathcal H(q,p,\alpha)}\leq \frac{2\delta}{3}+\sum_{n\in N(m_\delta)}\left|c_n\right|\left\|f_n\ast\varphi_{\epsilon\rho^{-1}_{n}}-f_n\right\|_{q',p'}.$$
By hypothesis, we have $1\leq p'\leq q'<\infty$ and therefore
$$\lim_{t\rightarrow 0}\left\|g\ast\varphi_t-g\right\|_{q',p'}=0,\ g\in(L^{q'},\ell^{p'})(\mathbb R^d).$$
So we have
$$\lim_{\epsilon\rightarrow 0}\left\|f_n\ast\varphi_{\epsilon\rho^{-1}_{n}}-f_n\right\|_{q',p'}=0,\ n\in N(m_\delta)$$
and therefore
$$\overline{\lim_{\epsilon\rightarrow 0}}\left\|f\ast\varphi_\epsilon-f\right\|_{\mathcal H(q,p,\alpha)}\leq\frac{2\delta}{3}.$$
This inequality being true for any real number $\delta>0$, we actually have
$\lim_{\epsilon\rightarrow 0}\left\|f\ast\varphi_\epsilon-f\right\|_{\mathcal H(q,p,\alpha)}=0$.

\item  Approximating a given function in
$ {\mathcal H(q,p,\alpha)}$ first by some function with compact support and then convolving it by some compactly supported, infinitely differentiable test function provides an approximation by a test function, which also belongs to the Schwartz space.
\item For 
 $1<q\leq p\leq\infty$,
we have $1\leq p'\leq q'<\infty$. But it is well known that $(L^{q'},\ell^{p'})$ is separable. Thus the result follows from the density of $(L^{q'},\ell^{p'})$  in $ {\mathcal H(q,p,\alpha)}$.
\end{enumerate}
\epf

For $1<q\leq\alpha\leq p\leq\infty$, we have defined a pre-dual of the Fofana space $(L^q,\ell^p)^\alpha$ using some atomic decomposition method developed by Feichtinger, and proved that for $p=\infty$ and $q < \alpha$, the pre-dual of classical Morrey space is embedded in our space. But our further goal is to find the dual space of $(L^q,\ell^p)^\alpha$ and their interpolation spaces. This appears not as an easy task
and has to be left for future work. 
\vspace{3mm}

\section{Conclusion}

In summary we have used techniques concerning minimal invariant Banach spaces of functions
in order to characterize the pre-dual of certain Fofana spaces which have not been known so far.
Starting from a characterization of a Fofana space as a (dense) subspace of a Wiener amalgam
space under a certain group of (suitably normalized) dilation operators one can generate
the predual space by starting from the predual of the mentioned Wiener amalgam spaces
and then describing the predual via atomic decompositions, using the (adjoint) group
of dilation operators.

{\it Acknowledgement:}  During the period of the preparation of the material for this manuscript (spring and summer of 2018) the first author held a guest position at TU Muenich, Dept. of Theoretical Information Sciences (H.~Boche). The second author is thankful to Fofana Ibrahim for drawing his attention to the separability of the predual, and many helpful discussions. 

\vspace{5mm}

\end{document}